\documentclass[11pt,twoside]{article}

\usepackage{epsfig,amsfonts,color}
\usepackage{amsmath, bbm}

\usepackage{amssymb, palatino, geometry,url}
\usepackage{algorithmic}
\usepackage[noresetcount,lined,boxed]{algorithm2e} 

\usepackage[colorlinks=true,linkcolor=blue,citecolor=blue,urlcolor=blue]{hyperref}
\usepackage{subcaption}
\usepackage{amsthm}
\newtheorem{thm}{Theorem}

\usepackage{fancyhdr}
\newcommand{\be}{\begin{equation}}
\newcommand{\ee}{\end{equation}}
\newcommand{\bea}{\begin{eqnarray}}
\newcommand{\eea}{\end{eqnarray}}

\newcommand{\bvec}{\left(\begin{array}{c}}
	\newcommand{\evec}{\end{array}\right)}
\newcommand{\bsub}{\begin{subequations}}
	\newcommand{\esub}{\end{subequations}}

\usepackage{lineno}
\usepackage{float}
\usepackage[section]{placeins}
\usepackage[table]{xcolor}

\title{Measuring and Optimizing System Reliability: \\A Stochastic Programming Approach}
\author{Joshua L. Pulsipher and Victor M. Zavala\thanks{Corresponding Author: victor.zavala@wisc.edu}\\
	{\small Department of Chemical and Biological Engineering}\\
	{\small \;University of Wisconsin-Madison, 1415 Engineering Dr, Madison, WI 53706, USA}}
\date{}

\begin{document}
	
\maketitle

\begin{abstract}
We propose a computational framework to quantify (measure) and to optimize the reliability of complex systems. The approach uses a graph representation of the system that is subject to random failures of its components (nodes and edges). Under this setting, reliability is defined as the probability of finding a path between sources and sink nodes under random component failures and we show that this measure can be computed by solving a stochastic mixed-integer program. The stochastic programming setting allows us to account for system constraints and general probability distributions to characterize failures and allows us to derive optimization formulations that identify designs of maximum reliability.  We also propose a strategy to approximately solve these problems in a scalable manner by using purely continuous formulations. 
\end{abstract}


\section{Introduction}

In this work, we investigate the problem of quantifying the reliability of complex systems and of designing systems of maximum reliability. Such problems have a wide range of applications such as supply chains, transportation networks, energy networks, process networks, sensor networks, and control networks  \cite{kim2013network}. In these applications, it is vital to design systems that maintain functionality in the face of natural and man-made events (e.g., mechanical failures, power outages, weather, and cyber-attacks) \cite{yan2012survey}.  Despite its practical importance, quantifying the reliability of complex systems remains a technical challenge. 

Reliability has been traditionally defined as the probability that a system remains functional under component failures  \cite{ogunnaike2009random}. The most prominent model used in industry to quantify reliability is based on so-called reliability block diagrams (RBDs). Here, the system is modeled as a network (a directed graph) of series/parallel paths in which each path has a single source and sink node. The system is said to function under a given failure if there exists at least one  path between the source and the sink node. The RBD approach exploits the simple topology of series/parallel systems to analytically compute the reliability of the overall system from the reliability of its individual components  \cite{thomaidis1994integration}. Here, it is also implicitly assumed that the probability of failure for every component can be chracterized using the same probability distribution. The availability of an analytical measure facilities the design of systems of maximum reliability \cite{ye2018mixed}. Unfortunately, the RBD approach is difficult to apply to more complex settings that involve, for instance, topologies with multiple source and sink nodes and loops and components with different probability distributions. As a result, analytical reliability measures cannot be easily derived under such settings.

The recursive decomposition algorithm (DFA) is a technique that aims to quantify reliability of more complex network topologies by systematically exploring paths between source and sink nodes  \cite{bistouni2014analyzing}.  This approach is more general but is not amenable for design tasks. Simulation-based approaches such as Monte Carlo (MC) sampling provide a general approach to quantify reliability. These approaches estimate reliability by ``probing" the system against failure scenarios and from this determine the probability that the system remains functional \cite{li2013reliability}. This approach is computationally more expensive than the analytical RBD approach because it requires repetitive simulations but can also enable the use of a wide range of stochastic programming formulations and solution techniques \cite{luedtke2008sample}. Specifically, we show that reliability can be computed by solving a stochastic mixed-integer program. This framework allows us to handle arbitrary system topologies, probability distributions to characterize different types of failures, and system constraints. Moreover, the stochastic program can be easily incorporated within optimal design formulations. We also provide evidence that accurate solutions for large systems can be obtained by solving purely continuous relaxations. 

The paper is structured as follows: Section \ref{sec:problem_def} establishes the definition of reliability guiding this work and introduces basic notation. Section \ref{sec:theory} provides stochastic programming formulations to compute reliability and to design systems with maximum reliability. Section \ref{sec:cases} presents case studies. Section \ref{sec:conclusion} provides concluding remarks.

\section{Problem Definition and Setting} \label{sec:problem_def}
In this section, we present a general graph abstraction to model complex systems. This abstraction is used to motivate and define reliability measures. 

\subsection{Graph Abstraction and Model}
\begin{figure}[!hbt]
	\centering
	\includegraphics[width=0.3\textwidth]{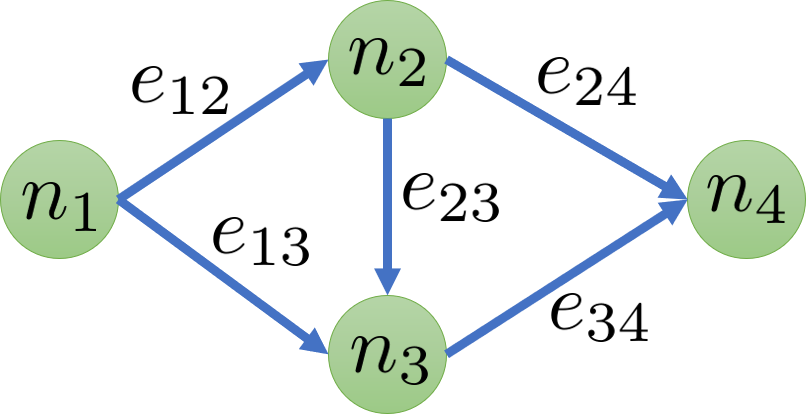}
	\caption{Representation of a system as a directed graph  with node set $\mathcal{N}=\{n_1,n_2,n_3,n_4\}$ and edge set $\mathcal{E}=\{e_{12},e_{13},e_{23},e_{24},e_{34}\}$.}
	\label{fig:graph_example}
\end{figure}

We model a system as a directed graph $\mathcal{G}(\mathcal{N}, \mathcal{E})$ with components $\mathcal{N}$ (nodes) and $\mathcal{E}$ (edges). We use $n \in \mathcal{N}$ and $e \in \mathcal{E}$ to represent specific nodes and edges in the graph, respectively. The set of edges originating at node $n$ is denoted as $\mathcal{E}_{in}(n) \subseteq \mathcal{E}$ and  the set of edges ending a node $n$ is denoted as $\mathcal{E}_{out}(n) \subseteq \mathcal{E}$. The set of supporting nodes for an edge $e$ (the pair of nodes connected by the edge) is denoted $\mathcal{N}(e) \subseteq \mathcal{N}$. A schematic representation of the graph notation is provided in Figure \ref{fig:graph_example}. 

The topology of the system $\mathcal{G}(\mathcal{N}, \mathcal{E})$ is encoded in the incidence matrix $A \in \mathbb{R}^{|\mathcal{N}| \times |\mathcal{E}|}$ where $A_{ne} = 1$ if $e \in \mathcal{E}_{in}(n)$, $A_{ne} = -1$ if $e \in \mathcal{E}_{out}(n)$, or $A_{ne} = 0$ otherwise. The nominal topology $A$ is subject to {\em failures} of its components (nodes and edges); as such, we define the {\em perturbed} incidence matrix as a random matrix $A(\xi_\mathcal{N}, \xi_\mathcal{E})$. Here, $\xi_\mathcal{N} \in \mathbb{R}^{|\mathcal{N}|}$ is the realization of a discrete (binary) random vector that indicates the set of nodes that function ($\xi_{\mathcal{N}, n} = 1$ if node $n$ functions) or do not function ($\xi_{\mathcal{N}, n} = 0$ if $n$ does not function). Similarly, $\xi_\mathcal{E} \in \mathbb{R}^{|\mathcal{E}|}$ denotes the realization of a binary random vector that indicates the set of nodes that function ($\xi_{\mathcal{E}, e} = 1$) or do not function ($\xi_{\mathcal{E}, e} = 0$). Under these definitions, the perturbed incidence matrix under realization $\xi:=(\xi_\mathcal{N}, \xi_\mathcal{E})$ can be computed as: 
\begin{equation}
A(\xi) := \Xi_\mathcal{N} {A} \Xi_\mathcal{E}
\label{eq:stochastic_incidence}
\end{equation}
where $\Xi_\mathcal{N} \in \mathbb{R}^{|\mathcal{N}| \times |\mathcal{N}|}$, $\Xi_\mathcal{E} \in \mathbb{R}^{|\mathcal{E}| \times |\mathcal{E}|}$ are diagonal matrices of the form $\Xi_\mathcal{N}=\textrm{diag}(\xi_{\mathcal{N}})$ and $\Xi_\mathcal{E}=\textrm{diag}(\xi_{\mathcal{E}})$, respectively. In a stochastic programming context, one can interpret $A(\xi)$ as a random technology matrix \cite{birge2011introduction}. The elements of the perturbed incidence matrix can also be written as:
\begin{align}
A_{ne}(\xi)= A_{ne} \cdot \xi_{\mathcal{N}, n} \cdot \xi_{\mathcal{E}, e},\; n\in\mathcal{N},\,e\in\mathcal{E}.
\end{align}
In other words, $A_{ne}(\xi)=0$ (entry does not exist) if either node $n$ or edge $e$ fails (do not exist) in scenario $\xi$. 

We use a network flow model to represent paths between nodes. Specifically, we define a set of {\em source} nodes as $\mathcal{N}_{so} \subseteq \mathcal{N}$ with associated source flows $d_n>0$, a set of {\em sink} nodes as $\mathcal{N}_{si}\subseteq \mathcal{N}$ with associated sink flows $d_n<0$, and a set of {\em relay} nodes as $\mathcal{N}_{re} \subseteq \mathcal{N}$ with associated flows $d_n=0$.  We observe that the source and sink flows are {\em fixed}. Under these definitions, the network flow representation can be expressed as:
\begin{align}
\sum_{e\in \mathcal{E}} A_{ne}(\xi) z_e + d_{n}  = 0,\; n \in \mathcal{N}
\end{align}
where $z_e\in \mathbb{R}_+$ is the flow along edge $e\in\mathcal{E}$. The network flow model can also be expressed in compact form as:
\begin{align}
A(\xi) z + d  = 0.
\end{align}
In our framework, we expand this basic network flow model to capture the possibility of readjusting flows in order to maintain system functionality. This can be done by allowing some nodes $\mathcal{N}_u\subseteq \mathcal{N}$ to have controllable flows $u_n\in \mathbb{R}_+$. Moreover, in many applications, the edge flows $z$ and the controls $u$ have physical meaning and are thus subject to constraints; we capture such constraints by using feasible sets $\mathcal{Z} \subseteq \mathbb{R}^{|\mathcal{E}|}$ and $\mathcal{U}\subseteq \mathbb{R}^{|\mathcal{N}_{u}|}$. With this, we define the extended network flow model as:
\begin{subequations}
\begin{align}
A(\xi) z + u + d = 0 \\ 
u \in \mathcal{U},\; z \in \mathcal{Z}.
\end{align}
\end{subequations}
In this representation, the set $\mathcal{U}$ is constructed in a way that it restricts control at certain nodes. For instance, we consider the box control set:
\begin{align}
\mathcal{U}=\{u\,:\,u_n=0,\; n\notin \mathcal{N}_u\; \&\; \underline{u}_n\leq u_n\leq \overline{u}_n,\; n\in \mathcal{N}_u\}. 
\end{align}
For simplicity, we assume that the feasible set for flows is also a box set of the form:
 \begin{align}
\mathcal{Z}=\{z\,:\underline{z}_e\leq z_e\leq \overline{z}_e,\; e\in \mathcal{E}\}. 
\end{align}

\subsection{Reliability Measures} \label{sec:trad_reliability}
A {\em reliability measure} seeks to quantify the probability that a system remains functional under random component failures. Under a graph representation, the system is said to be functional if there exists at least one path that connects each sink node to a source node. For a particular realization $\xi$ (with associated topology $A(\xi)$) and in the absence of controls and constraints, the functionality of a system can be checked by using the {\em reliability function}:
\begin{equation}
\psi(A,\xi) := \begin{cases} 1 &\text{if} \ \exists\, z: A(\xi) z + d  = 0 \\ 0 &\text{otherwise}. \end{cases}
\label{eq:trad_feasibility}
\end{equation} 
This function uses the network flow representation to check if there exist a set of flows $z$ that connect sinks and sources. This is based on the observation that, if a path does not exist between a sink and at least one source node (e.g., the network becomes disconnected in a given failure scenario), then there is no set of flows $z$ that satisfies the flow constraint $A(\xi) z+ d=0$. 

The traditional definition of reliability does not account for constraints  and does not account for the possibility to control flows. To account for these features, we extend the {\em reliability function} as:
\begin{equation}
\psi(A,\xi,\mathcal{Z},\mathcal{U}) := \begin{cases} 1 &\text{if} \ \exists\, z\in \mathcal{Z},u\in\mathcal{U}: A(\xi) z + u +d  = 0 \\ 0 &\text{otherwise}. \end{cases}
\label{eq:full_feasibility}
\end{equation} 
We use this extended function to define the {\em reliability measure}:
\begin{equation}
R(A,\mathcal{Z},\mathcal{U}) := \mathbb{P}(\psi({A}(\xi),\mathcal{Z},\mathcal{U}) = 1).
\label{eq:trad_reliability}
\end{equation} 
This measure function is the probability that the system remains functional.  A similar measure has been proposed to measure system {\em flexibility} which, in our setting, would represent the ability of a system to withstand perturbations in the source and sink flows $d$ (exogenous disturbances) \cite{straub1993design,bansal1998flexibility,swaney1985index,pulsipher2018mixed}. Therefore, we highlight that a key distinction between flexibility and reliability is that the former deals with continuous perturbations while the later deals with discrete perturbations. 

\subsection{Designs of Maximum Reliability}
We are interested in using the reliability measure to find system designs that maximize reliability. In this task, one often needs to trade-off cost $c(A,\mathcal{Z},\mathcal{U})$ and reliability, giving rise to the abstract problem:
\begin{equation}
\begin{aligned}
&&\max_{{A}, \mathcal{Z}, \mathcal{U}} &&&R(A, \mathcal{Z}, \mathcal{U}) \\ 
&&\text{s.t.} &&& c(A, \mathcal{Z}, \mathcal{U}) \leq \epsilon
\end{aligned}
\label{eq:epsilon_design}
\end{equation} 
where $\epsilon\in\mathbb{R}$ is a cost budget that is spanned to find Pareto pairs $(c^*, R^*)$.  We highlight the dependence of the cost measure and reliability measure on the topological design (given by the incidence matrix $A$) and on the operational design (given by the constraint sets $\mathcal{Z},\mathcal{U}$).

\section{Stochastic Programming Formulations} \label{sec:theory}
In this section, we provide stochastic programming formulations to compute the proposed reliability measure and to design systems of maximum reliability. We show that these formulations can be easily derived from the network flow representation of the system. 

\subsection{Computing the Reliability Measure}
We motivate the discussion by considering a simple setting with a single-input and single-output graph. Under this setting, the sets $\mathcal{N}_{so}$ and $\mathcal{N}_{si}$ are singletons and thus the system is said to remain functional if there exists at least one path between the source and the sink node. Equivalently, given a fixed source flow, the system is functional if we can find a set of edge flows that satisfy the fixed sink flow. Under this logic, we can compute $\psi(A,\xi)$ by finding a feasible solution for a network flow problem and this problem can be cast as a mixed-integer linear program (MILP) of the form:
\begin{equation}
\begin{aligned}
\psi(A,\xi)& =&\max_{y, z} &&&  (1-y) \\
&&\text{s.t.} &&& \sum_{e\in \mathcal{E}} A_{ne}(\xi) z_e = 0, && n \in \mathcal{N}_{re}\\
&&&&& \sum_{e\in \mathcal{E}} A_{ne}(\xi) z_e + d_n \cdot (1-y)=0, && n \in \mathcal{N}_{so}\\
&&&&& \sum_{e\in \mathcal{E}} A_{ne}(\xi) z_e +d_n \cdot(1-y) = 0, && n \in \mathcal{N}_{si}\\
&&&&& z_e \geq 0, && e\in \mathcal{E} \\
&&&&& y \in \{0, 1\}.
\end{aligned}
\label{eq:single_feasibility}
\end{equation}
Here, we arbitrarily set the source and sink flows to $d_n=1$ and $d_n=-1$, respectively. This is done without loss of generality because the flows do not necessarily have physical meaning (in more general settings they might have meaning). We use the binary variable $y\in \{0,1\}$ to relax the balances at the source and sink nodes (i.e., if $y=0$ then the network flow system has a feasible solution and if $y=1$ then it does not). If the network flow system does not have a solution then we obtain the trivial flow solution $z_e=0$ for all $e\in\mathcal{E}$. We thus have that the reliability measure is given by $\psi(A,\xi)=1-y^*$ and we note that the maximization problem is equivalent to minimize $y$. The MILP can be relaxed by setting $0 \leq y \leq 1$; interestingly, this LP is guaranteed to deliver an optimal (binary) solution for the MILP (see Appendix). 

Problem \eqref{eq:single_feasibility} can be easily generalized to compute the reliability measure for graphs with multiple sources and sinks and with controllable flows. This can be done by solving the MILP:
\begin{equation}
\begin{aligned}
\psi(A,\xi,\mathcal{Z},\mathcal{U})& =&\max_{y, z, u} &&& (1 - y) \\
&&\text{s.t.} &&& \sum_{e\in \mathcal{E}} A_{ne}(\xi) z_e = 0, && n \in \mathcal{N}_r\\
&&&&& \sum_{e\in \mathcal{E}} A_{ne}(\xi) z_e + d_n \cdot(1-y)+u_{n} = 0, && n \in \mathcal{N}_{so}\\
&&&&& \sum_{e\in \mathcal{E}} A_{ne}(\xi) z_e +d_n\cdot (1-y)+ u_{n} = 0, && n \in \mathcal{N}_{si}\\
&&&&&  u\in\mathcal{U}, \; z\in\mathcal{Z}\\
&&&&& y \in \{0,1\}.
\end{aligned}
\label{eq:multiple_feasibility}
\end{equation}
This MILP determines if all the sink flows can be satisfied via the source flows (i.e., each sink has at least one path to a source); this is true whenever $y = 0$ (which indicates that none of the source and/or sink nodes needs to be relaxed to achieve a feasible solution). 

The MILP representation of the reliability function reveals that the measure $R(A, \mathcal{Z},\mathcal{U})$ is a {\em joint chance constraint}. This chance constraint can be approximated using MC samples $\xi^k, \ k \in \mathcal{K}$ as \cite{kim2015guide}:
\begin{equation}
R(A,\mathcal{Z},\mathcal{U}) \approx \frac{1}{|\mathcal{K}|} \sum_{k \in \mathcal{K}} \psi(A,\xi^k,\mathcal{Z},\mathcal{U}).
\label{eq:r_definition}
\end{equation}
By the law of large numbers, this sample average approximation becomes asymptotically exact as the number of samples increases  \cite{hsu1947complete}; moreover, the approximation converges exponentially \cite{kleywegt2002sample}. Combining problems \eqref{eq:r_definition} and \eqref{eq:multiple_feasibility}, we obtain the following  approximation of the reliability measure:
\begin{equation}
\begin{aligned}
R(A,\mathcal{Z},\mathcal{U}) & \approx&\max_{y^k, z^k, u^k} &&& \frac{1}{|K|} \sum_{k \in \mathcal{K}} (1 - y^k) \\
&&\text{s.t.} &&& \sum_{e\in \mathcal{E}} A_{ne}(\xi^k)z_e^k = 0, && n \in \mathcal{N}_{re}, \ k \in \mathcal{K}\\
&&&&& \sum_{e\in \mathcal{E}} A_{ne}(\xi^k) z_e^k + d_n\cdot(1-y^k)+u_n^k = 0, && n \in \mathcal{N}_{so}, \ k \in \mathcal{K}\\
&&&&& \sum_{e\in \mathcal{E}} A_{ne}(\xi^k) z_e^k +d_n\cdot(1-y^k)+u_n^k= 0 , && n \in \mathcal{N}_{si}, \ k \in \mathcal{K}\\
&&&&& z^k\in\mathcal{Z},\; u^k\in \mathcal{U}, && k \in \mathcal{K}\\ 
&&&&& y^k \in \{0, 1\}, && k \in \mathcal{K}.
\end{aligned}
\label{eq:r_formulation}
\end{equation}
This problem is {\em fully decoupled} in the MC samples $k \in \mathcal{K}$ and thus can be trivially parallelized. It has been recently reported that a continuous relaxation of this problem (in combination with an appropriate rounding strategy)  provides high quality approximations of the exact solution \cite{pulsipher2019scalable}. Specifically, we can relax $y^k\in \{0,1\}$ to $0 \leq y^k \leq 1$ and then round the optimized relaxed $y^{k*}$ values to 1 if they are nonzero. This approach is analogous to employing slack variables to identify active and inactive sets of constraints.  In the following section we provide numerical evidence that this relaxation approach is effective. The exact relaxation result for the simple reliability problem \eqref{eq:single_feasibility} provides some intuition as to why this happens. However, establishing a theoretical justification in a more complex setting with constraints and controllable flows is difficult and is left as a topic of future work.

The MILP representation can be extended in a number of ways to capture desirable decision-making logic. For instance, one might want to relax the requirement that paths must exist to all sink nodes and instead require that only a subset of nodes are reachable. This can be done by introducing binary variables for all sink nodes $y_n^k$ and by solving the problem:
\begin{equation}
\begin{aligned}
R(A,\mathcal{Z},\mathcal{U}) & \approx&\max_{y^k, z^k, u^k} &&& \frac{1}{|K|} \sum_{k \in \mathcal{K}} L(y^k) \\
&&\text{s.t.} &&& \sum_{e\in \mathcal{E}} A_{ne}(\xi^k)z_e^k = 0, && n \in \mathcal{N}_{re}, \ k \in \mathcal{K}\\
&&&&& \sum_{e\in \mathcal{E}} A_{ne}(\xi^k) z_e^k + d_n\cdot(1-y_n^k)+u_n^k = 0, && n \in \mathcal{N}_{so}, \ k \in \mathcal{K}\\
&&&&& \sum_{e\in \mathcal{E}} A_{ne}(\xi^k) z_e^k +d_n\cdot(1-y_n^k)+u_n^k=0 , && n \in \mathcal{N}_{si}, \ k \in \mathcal{K}\\
&&&&& z^k\in\mathcal{Z},\; u^k\in \mathcal{U}, && k \in \mathcal{K}\\ 
&&&&& y^k_n \in \{0, 1\}, && k \in \mathcal{K},\; n\in\mathcal{N}_{si} \cup \mathcal{N}_{so}.
\end{aligned}
\label{eq:r_formulation2}
\end{equation}
Here, $L(y^k)$ is a logic function which is set to one if a subset of sinks of interest are reachable (or is set to zero otherwise).  

\subsection{Optimal Design} \label{sec:design}
The design problem \eqref{eq:epsilon_design} aims to make topological and capacity changes to a nominal network in order to maximize reliability (under a given cost budget).  To formulate this problem, we recall that the base topology of the system is given by the graph $\mathcal{G}(\mathcal{N},\mathcal{E})$ with associated incidence matrix $A$, nodes $\mathcal{N}$, and $\mathcal{E}$. Our goal is this formulation to expand the number of edges in order to maximize reliability. This is done by defining an expanded set of edges $\bar{\mathcal{E}}$ such that $\mathcal{E}\subset\bar{\mathcal{E}}$.  The expanded set of edges has an associated incidence matrix $\bar{A}$. In other words, the new incidence matrix has an expanded set of connections between the nodes. We represent the added set of edges as $\hat{\mathcal{E}}:=\bar{\mathcal{E}}\setminus\mathcal{E}$. In our design problem, we also seek to expand the set of feasible edge flows and control flows (to model capacity expansions). The design problem is cast as the following MILP: 
\begin{equation}
\begin{aligned}
&&\max_{v,\underline{z}, \overline{z}, \underline{u}, \overline{u}, z^k,y^k, u^k} &&& \frac{1}{|K|} \sum_{k \in \mathcal{K}} L(y^k) \\ 
&&\text{s.t.} &&& c(v,\underline{z}, \overline{z}, \underline{u}, \overline{u}) \leq \epsilon \\
&&&&& \sum_{e\in \bar{\mathcal{E}}} A_{ne}^kz_e^k = 0, && n \in \mathcal{N}_{re}, \ k \in \mathcal{K}\\
&&&&& \sum_{e\in \bar{\mathcal{E}}} A_{ne}^k z_e^k + d_n\cdot(1-y_n^k)+u_n^k = 0, && n \in \mathcal{N}_{so}, \ k \in \mathcal{K}\\
&&&&& \sum_{e\in \bar{\mathcal{E}}} A_{ne}^k z_e^k +d_n\cdot(1-y_n^k)+u_n^k=0 , && n \in \mathcal{N}_{si}, \ k \in \mathcal{K}\\
&&&&& A_{ne}^k = \bar{A}_{ne}\cdot \xi_{\mathcal{N}, n}^k\cdot \xi_{\mathcal{E}, e}^k, && n \in \mathcal{N}, \ e \in \bar{\mathcal{E}}, \ k \in \mathcal{K} \\
&&&&& A_{ne}^k = \bar{A}_{ne}\cdot \xi_{\mathcal{N}, n}^k\cdot \xi_{\mathcal{E}, e}^k \cdot v_{e}, && n \in \mathcal{N}, \ e \in \hat{\mathcal{E}}, \ k \in \mathcal{K} \\
&&&&& \underline{z} \leq z^k \leq \overline{z},\;  \underline{u} \leq u^k \leq \overline{u}, &&  k \in \mathcal{K} \\
&&&&& y^k_n \in \{0, 1\},&& n\in\mathcal{N}_{si}\cup\mathcal{N}_{so},\; k \in \mathcal{K} \\
&&&&& v_e \in \{0, 1\},&& e\in \hat{\mathcal{E}} \\
&&&&& \underline{z},\overline{z} \in \overline{\mathcal{Z}},\; \underline{u},\overline{u} \in \overline{\mathcal{U}}.
\end{aligned}
\label{eq:mip_epsilon_design1}
\end{equation} 
Here, the sets $\overline{\mathcal{Z}}$ and $\overline{\mathcal{U}}$ include possible design values for flow and control bounds. Also, $v \in \{0, 1\}^{|\hat{\mathcal{E}}|}$ denote topological design variables for selecting which of the candidate edges are included in the new design (if none are added then $v_e=0$ for all $e\in\hat{\mathcal{E}}$ and the network retains its nominal topology). We note that the abstract design cost function $c(A,\mathcal{Z},\mathcal{U}$) can now be expressed in the parametric form $c(v,\underline{z}, \overline{z}, \underline{u}, \overline{u})$. The proposed design formulation seeks to highlight the modeling flexibility provided by the proposed stochastic programming framework.  

\section{Case Studies} \label{sec:cases}

We analyze the behavior of the proposed framework by applying it to distribution networks. We consider a simple three-node network and the IEEE-14 power distribution network. We also consider a simple parallel-series RBD system to illustrate how the proposed stochastic programming framework is consistent with the analytical RBD solution. All formulations are implemented in {\tt JuMP} 0.18.5 \cite{dunning2017jump} and are solved using {\tt Gurobi} 7.5.1  on a Intel\textregistered \, Core\texttrademark \, i7-7500U machine running at 2.90 GHz with 4 hardware threads and 16 GB of RAM. All results can be reproduced using the scripts provided in \url{https://github.com/zavalab/JuliaBox/tree/master/ReliableDesign}. 

\subsection{Reliability of Parallel-Series Systems}
\begin{figure}[!hbt]
	\centering
	\includegraphics[width=0.5\textwidth]{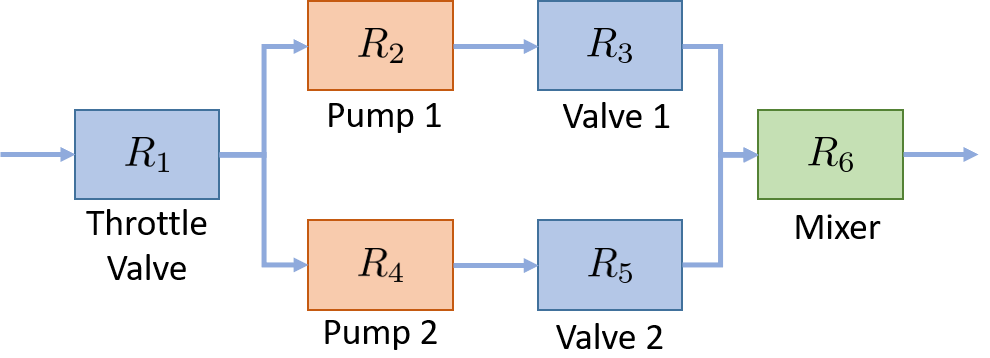}
	\caption{Reliability block diagram for a pump system.}
	\label{fig:rbd}
\end{figure}

We consider a simple parallel-series system to highlight that the stochastic programming approach is consistent. 
The system of interest is represented by the reliability block diagram shown in Figure \ref{fig:rbd}. This system seeks to pump a flow stream using two pumps and valves in parallel. The parallel design topology enhances the reliability of the system (compared to a topology with a single pump and valve). 

Traditional RBD methods can be leveraged to obtain an analytic representation for the overall reliability of the system since this system features a single source and sink \cite{bistouni2014analyzing}. In particular, the analytic reliability measure is computed by aggregating the component reliabilites according to their respective connectivities. Specifically, the reliability of $m$ components in series configurations are given by:
\begin{equation}
	R_s = \prod_i^m R_i.
\end{equation}
The reliability of a parallel configuration is given by:
\begin{equation}
	R_p = 1 - \prod_i^m (1 - R_i).
\end{equation}
Following these basic rules, the reliability of the system of interest can be computed as:
\begin{equation}
	R_{overall} = R_1 (1-(1-R_2R_3)(1-R_4R_5))R_6.
	\label{eq:rbd_analytic}
\end{equation}
For simplicity, we let each component lifetime be described by an exponential distribution with a mean lifetime of 100 years and we evaluate reliability after 5 years of operation. The reliability for each component is given by the exponential cumulative distribution function evaluated at 5 years (i.e., $R_i = \exp{(-5/100)}$). From Equation \eqref{eq:rbd_analytic} we thus obtain the overall reliability $R_{overall} = 89.66\%$. 

To demonstrate the equivalence of the proposed stochastic programming setting, we use MC samples drawn from the component distribution functions (a component fails if the lifetime is above the desired threshold of 5 years). We use a total of 10,000 MC samples and solve Problem \eqref{eq:r_formulation}. Here, we let the throttle valve be the source node and the mixer be the sink node with $d_n = 1$ and $d_n = -1$, respectively. Furthermore, we set $\mathcal{U} = \emptyset$ (no controls) and $\mathcal{Z} = \mathbb{R}_+^{|\mathcal{E}|}$. Using this approach, the reliability measure is $R(A, \mathcal{Z}, \mathcal{U}) = 89.69\%$, which is close to the analytical solution. 

\subsection{Network Models}
\begin{figure}[!hbt]
	\centering
	\includegraphics[width=0.4\textwidth]{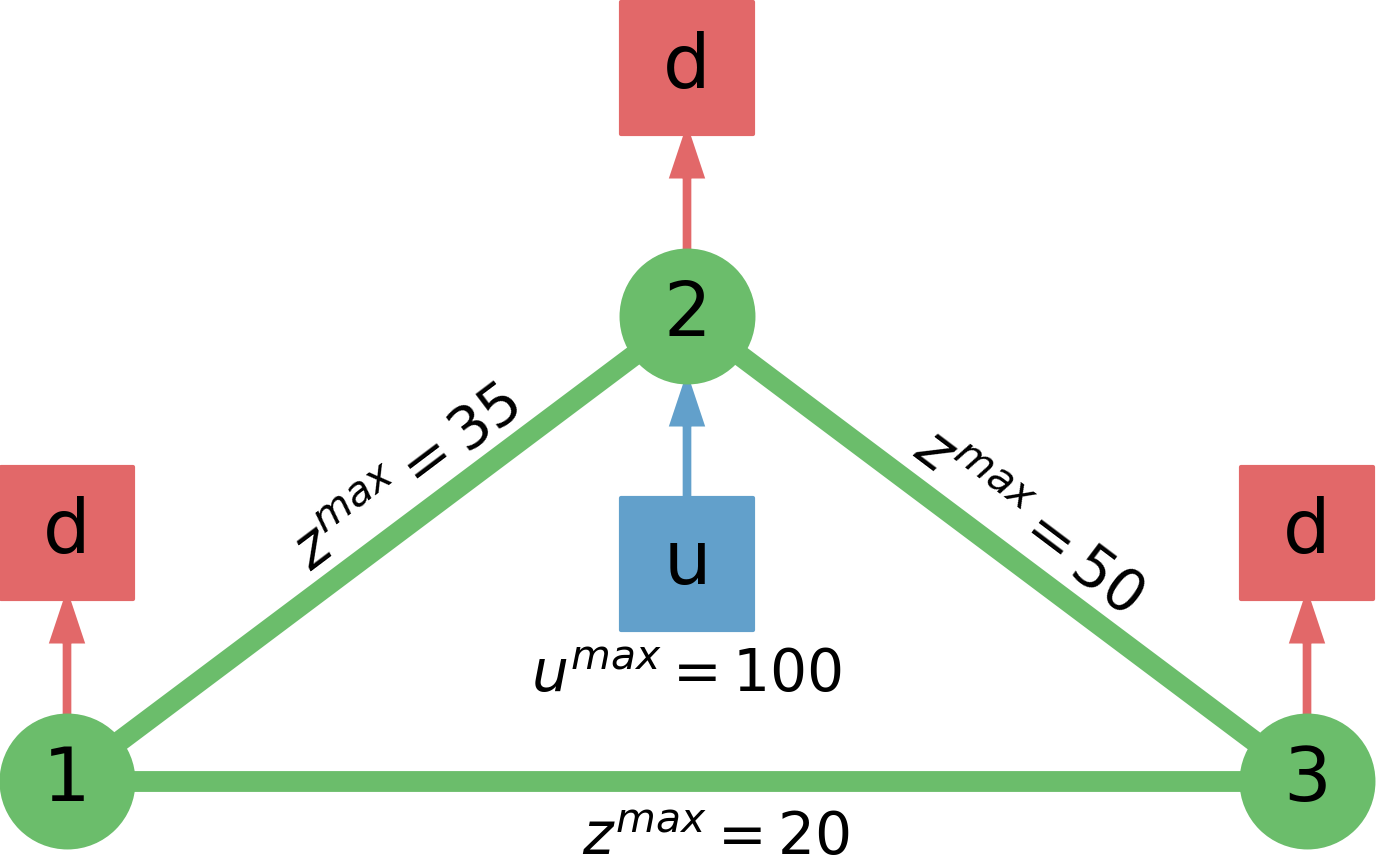}
	\caption{Schematic of 3-node distribution network.}
	\label{fig:3node}
\end{figure}

\begin{figure}[!hbt]
	\centering
	\includegraphics[width=0.6\textwidth]{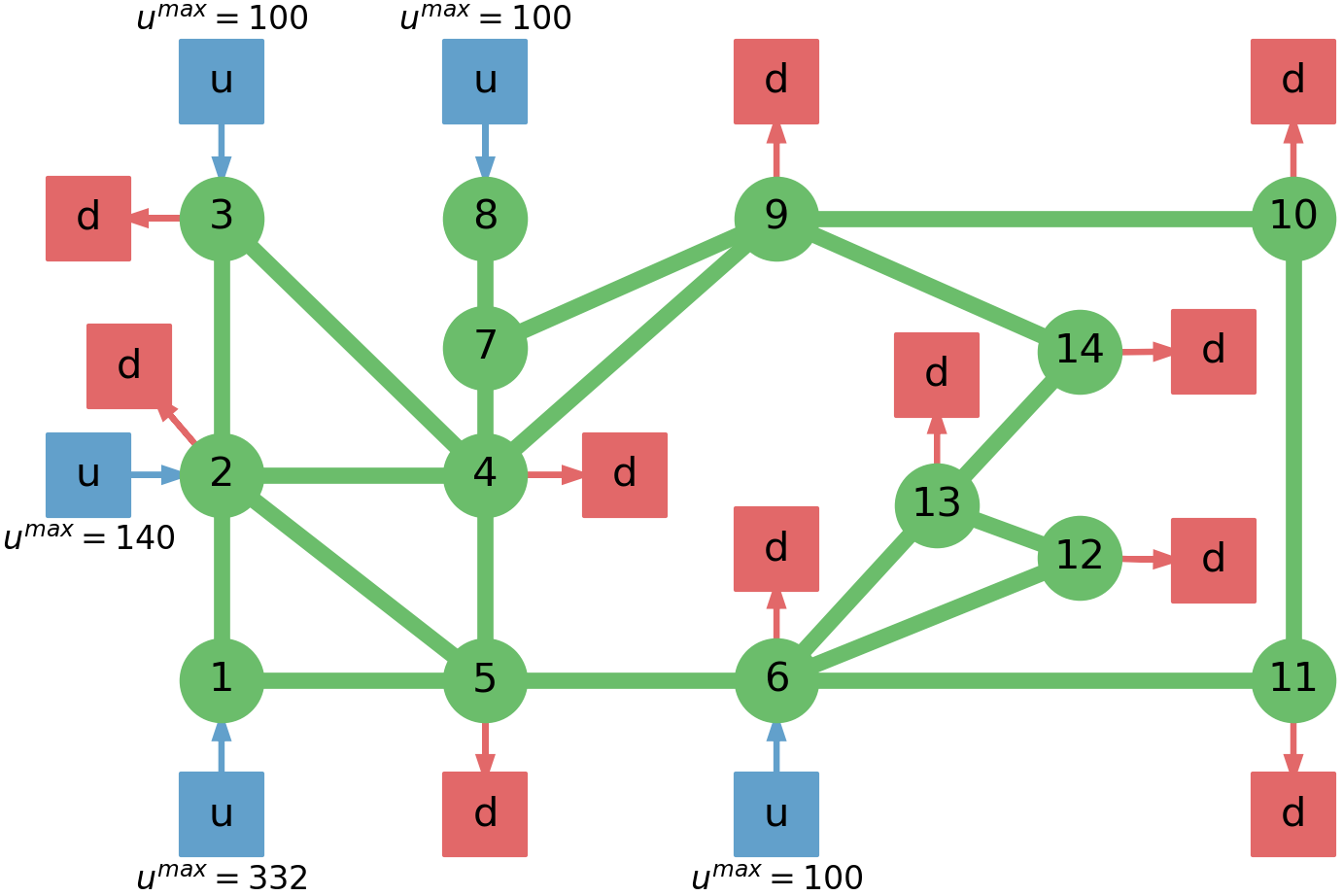}
	\caption{Schematic of IEEE 14-node distribution network.}
	\label{fig:ieee14}
\end{figure}

The systems under study are illustrated in Figures \ref{fig:3node} and \ref{fig:ieee14}. We consider a simple 3-node distribution network and the IEEE 14-node power network benchmark problem. In these cases, the sink nodes $n\in \mathcal{N}_{si}$ have a fixed flow $d_{n}$ and the source nodes $n\in\mathcal{N}_{so}$ are controllable with capacity $\bar{u}$. Furthermore, the edges have a finite capacity $\bar{z}$.  The 3-node network features a single source (a power plant) and three sink nodes (power consumers). The IEEE 14-node network exhibits a more complex topology with multiple sinks and sources. The data for this problem is obtained from MATPOWER \cite{zimmerman2010matpower}. 

For our design studies, we consider a cost function of the form:
\begin{equation}
c(v, \overline{z}, \overline{u}) := \sum_{e \in \mathcal{E}} (100\cdot v_e + \overline{z}_e) + \sum_{n \in \mathcal{N}_u} \overline{u}_{n}
\label{eq:cost}
\end{equation}
We consider random failures for relay nodes, source nodes, and sink nodes. Here, we model the lifetimes of the components as exponential random variables with mean lifetimes of 100, 80, and 50 years, respectively. MC failure scenarios $\xi_\mathcal{N}^k,  \ \xi_\mathcal{E}^k$ are generated by sampling the exponential distributions and the components are set to failure mode if their lifetime is above a certain threshold value. The thresholds for the 3-node and IEEE 14-node networks were set to 5 and 2 years, respectively.  

\subsection{Design for Maximum Reliability (Capacity Expansion)} \label{sec:mi_cap_design}

\begin{figure}[!hbt]
	\centering
	\includegraphics[width=0.6\textwidth]{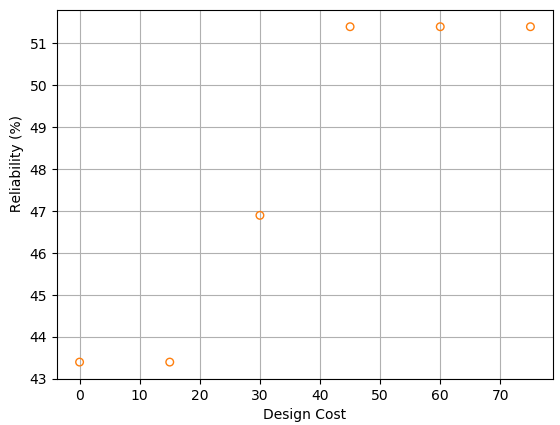}
	\caption{Pareto frontier for optimal capacity design of 3-node network using 1,000 MC samples.}
	\label{fig:3node_cap_mip}
\end{figure}

We first consider a design problem in which capacity is expanded (i.e., topological expansion variables $v$ are omitted). For the 3-node power network, we solve the MILP formulation to obtain 6 Pareto pairs and we use 1,000 MC samples. The Pareto pairs are plotted in Figure \ref{fig:3node_cap_mip}.  Here, we note that the Pareto frontier shows abrupt changes; this is because this system is simple and thus the solution space is small. A manifestation of this limited spaces is that the maximum possible reliability for this system is just 51.4\%. This indicates that, regardless of how much capacity is provided (unlimited budget), this network will never achieve a higher reliability because of its limiting topology. In other words, the only way to increase reliability is to add edges. 

\begin{figure}[!hbt]
	\centering
	\begin{subfigure}[t]{.4\textwidth}
		\centering
		\includegraphics[width=\textwidth]{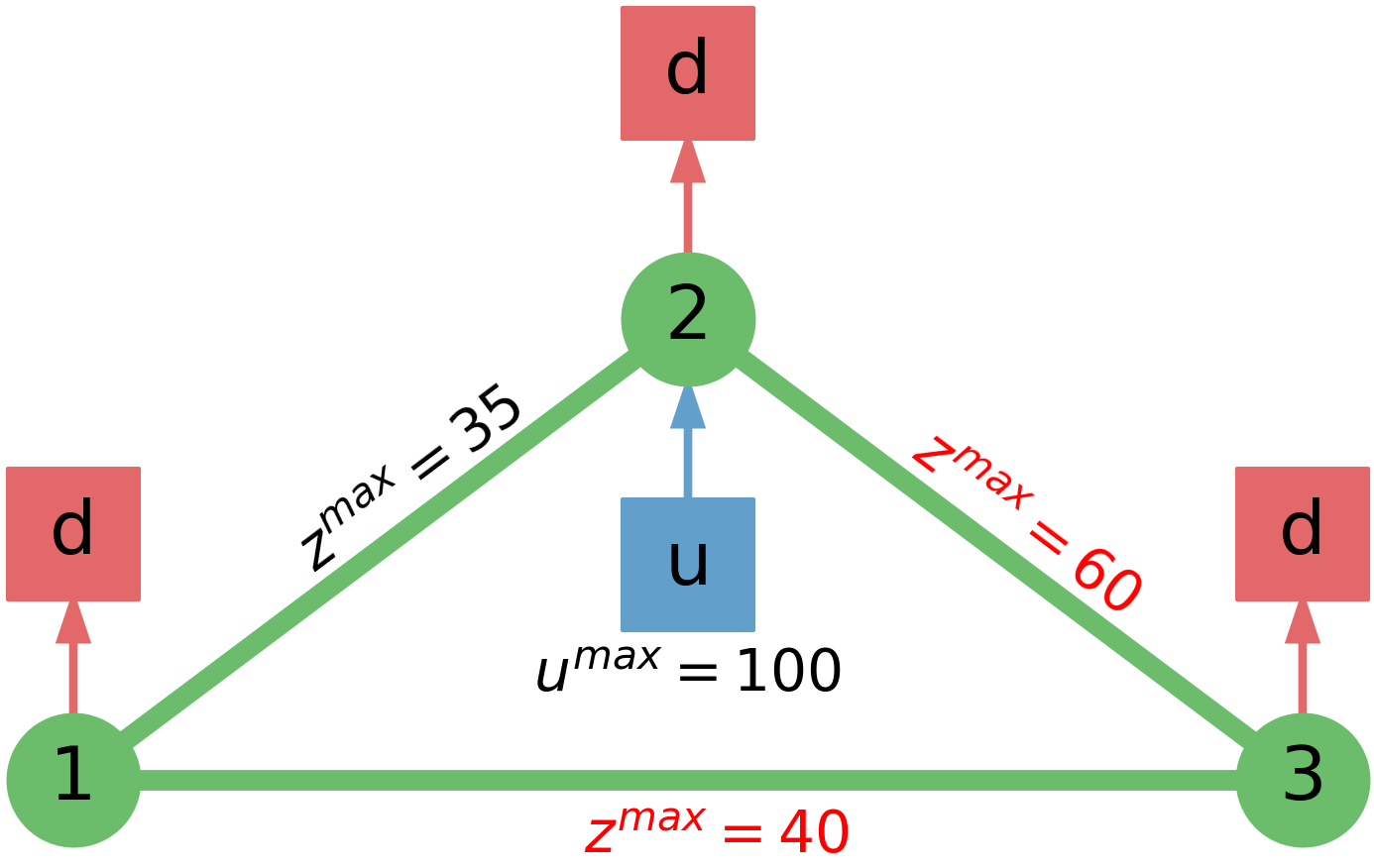}
	\end{subfigure}
	\caption{Schematic of the 3-node network corresponding to Pareto pair shown in Figure \ref{fig:3node_cap_mip} with a design cost of 30.}
	\label{fig:3node_cap_combined}
\end{figure}

We explore the Pareto solutions obtained with $\epsilon=30$ and $\epsilon=45$. The optimized capacities for these solutions are shown in Figure \ref{fig:3node_cap_combined}. Here, the increased capacities relative to the base design are highlighted in red. In the first case, enough capacity is added to the edges connecting nodes 2, 3, and 1 to permit the network to function in the event that the edges connecting nodes 1 and 2 fails. In the other case, enough capacity is added to the edges to permit feasible operation if any one edge fails.

\begin{figure}[!hbt]
	\centering
	\includegraphics[width=0.6\textwidth]{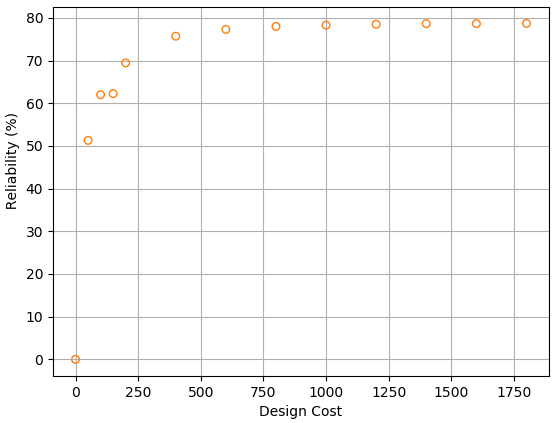}
	\caption{Pareto frontier for optimal design problem of IEEE 14-node network.}
	\label{fig:14node_cap_mip}
\end{figure}

We apply the same design formulation to the IEEE-14 power network problem. We compute a total of 13 Pareto pairs by varying the budget $\epsilon$ from 0 to 1800 and we use 2,000 MC samples. The solutions obtained with the MILP formulation are presented in Figure \ref{fig:14node_cap_mip}. We see that this system shows a smoother Pareto frontier because the solution space for this more complex system is larger. For this system the largest possible reliability is 78.7\% (this system has more degrees of freedom).

\begin{figure}[!hbt]
	\centering
	\includegraphics[width=0.6\textwidth]{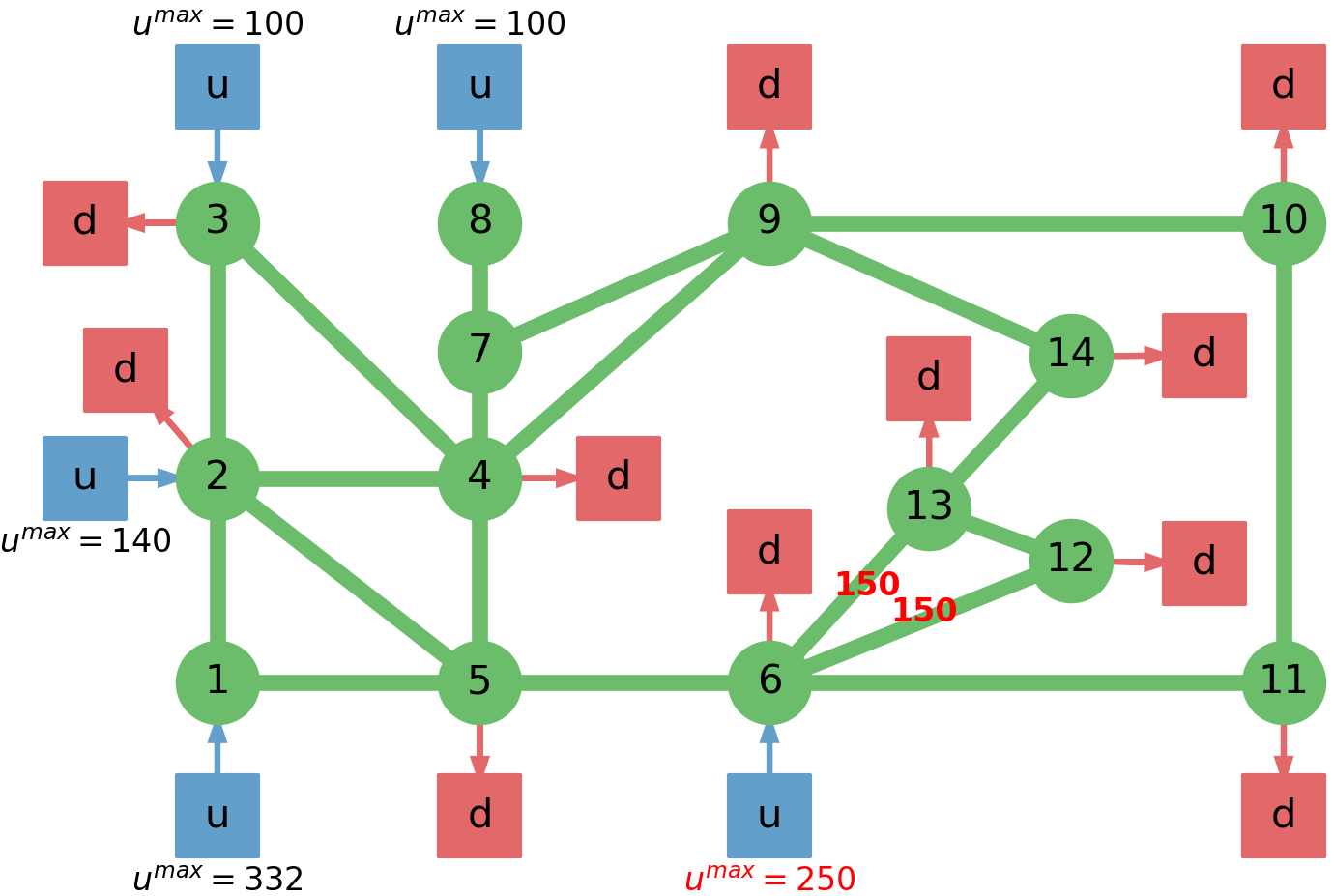}
	\caption{Schematic of the IEEE 14-node network corresponding to Pareto pair shown in Figure \ref{fig:14node_cap_mip} with a cost of 400.}
	\label{fig:14node_cap_design}
\end{figure}

The design obtained with a budget of $\epsilon=400$ is shown in Figure \ref{fig:14node_cap_design}. The expanded  capacities are highlighted in red. The capacity of the supplier attached to relay node 6 is significantly expanded (which occurs because it is the only supplier that serves the right side of the network). The capacities of two edges attached to node 6 are also increased such that the internal demands can be satisfied if either edge fails. It is interesting to note that these 3 simple changes to the network design significantly increase the overall reliability of the system (they increase it by 24.4\%). 

\subsection{Continuous Approximation for Design Problem} \label{sec:cont_cap_design}
\begin{figure}[!hbt]
	\centering
	\includegraphics[width=0.6\textwidth]{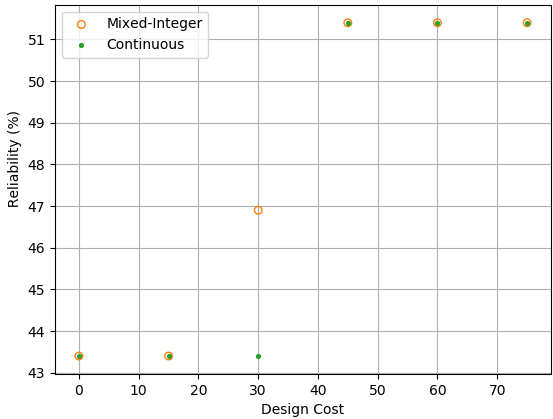}
	\caption{Pareto frontier for optimal capacity design of the 3-node network juxtaposing the pairs obtained from the full MILP formulation and its continuous relaxation.}
	\label{fig:3node_cap_cont}
\end{figure}

We consider a continuous relaxation of the design problem; here, integer solutions are obtained by solving the relaxation and then using simple rounding. This technique is first applied to the 3-node power network using the same samples and $\epsilon$ values considered above in Section \ref{sec:mi_cap_design}. In Figure \ref{fig:3node_cap_cont}, we juxtapose the resulting Pareto pairs. We observe that 5 out of 6 pairs are exactly recovered and a pair is underestimated. In \cite{pulsipher2019scalable} it is hypothesized that the quality of the approximations is the result of degeneracy associated with the joint chance constraint (i.e., multiple solutions yield the same optimal value). This simple network exhibits little degeneracy at that solution because its solution space is small. 

Table \ref{tab:3node_cap_data}  summarizes the performance of the MILP formulation and the continuous relaxation for the 3-node network. We observe that 5 of the 6 pairs are exactly equivalent since they have no differences in the active constraints. Also, the third pair only differs by 3.5\% (which is a small gap). For this small network, the solution times are negligible so the benefits of the relaxation are not obvious. 

\begin{table}[!htb]
	\caption{Performance results obtained for 3-node network using the MILP design formulation and continuous relaxation.}
	\begin{center}
		\begin{tabular}{|cccccc|}
			\hline	
		Cost & MILP R (\%) & LP R (\%) & MILP Time (s) & LP Time (s) &  Active Difference (\%) \\ \hline \hline 
			0           & 43.4         & 43.4       & 0.0467        & 0.0349      & 0   \\ 
			15          & 43.4         & 43.4       & 0.0400        & 0.0468      & 0   \\ 
			30          & 46.9         & 43.4       & 0.0407        & 0.0478      & 3.5 \\ 
			45          & 51.4         & 51.4       & 0.0562        & 0.0478      & 0   \\ 
			60          & 51.4         & 51.4       & 0.0526        & 0.0478      & 0   \\ 
			75          & 51.4         & 51.4       & 0.0443        & 0.0458      & 0   \\ 
			\hline
		\end{tabular}
	\end{center}
	\label{tab:3node_cap_data}
\end{table}

\begin{figure}[!hbt]
	\centering
	\includegraphics[width=0.6\textwidth]{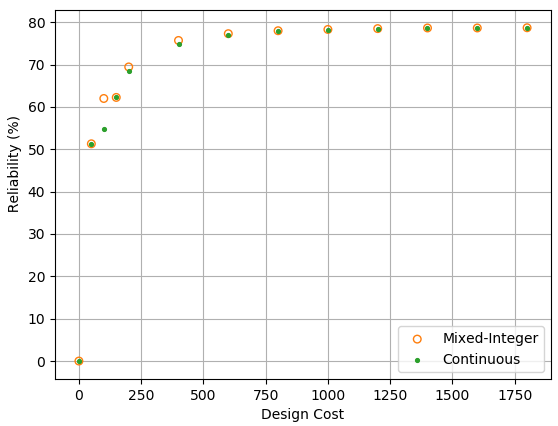}
	\caption{Pareto frontier for optimal capacity design of the IEEE 14-node network juxtaposing the pairs obtained from the full MILP formulation and its continuous relaxation.}
	\label{fig:14node_cap_cont}
\end{figure}

The relaxation strategy was also applied to the IEEE 14-node power network using the same conditions specified above in Section \ref{sec:mi_cap_design}. A juxtaposition of the Pareto pairs is shown in Figure \ref{fig:14node_cap_cont}. We observe that the frontier is approximated well (the majority of the pairs being exactly reproduced). Some of the minor discrepancies are attributed to numerical precision. A summary of the results is shown in Table \ref{tab:14node_cap_data}. With the exception of the third pair, the Pareto pairs only exhibit differences in the active constraints of less than 1\%. We can thus see that the relaxation indeed delivers high quality approximation. Importantly, we observe that the computational time is reduced by 96\%. This enables us to handle much larger networks than would be possible using the full MILP formulation. 

\begin{table}[!htb]
	\caption{The performance results obtained for the IEEE 14-node network using the mixed-integer and continuous capacity design formulations.}
	\begin{center}
		\begin{tabular}{|cccccc|}
			\hline	
			Design Cost (-) & MILP R (\%) & LP R (\%) & MILP Time (s) & LP Time (s) & Active Difference (\%) \\ \hline \hline 
			0     & 0     & 0     & 5.42  & 1.23  & 0 \\
			50    & 51.3  & 51.3  & 58.84 & 12.84 & 0 \\
			100   & 62    & 54.9  & 78.62 & 3.68  & 7.1 \\
			150   & 62.25 & 62.25 & 95.09 & 5.75  & 0.5 \\
			200   & 69.45 & 68.4  & 83.35 & 4.75  & 1.05 \\
			400   & 75.7  & 74.85 & 124.46 & 4.84  & 0.85 \\
			600   & 77.3  & 77.1  & 168.61 & 7.39  & 0.8 \\
			800   & 78    & 77.9  & 185.11 & 3.97  & 0.2 \\
			1000  & 78.3  & 78.1  & 205.75 & 7.2   & 0.5 \\
			1200  & 78.5  & 78.45 & 186.34 & 3.27  & 0.05 \\
			1400  & 78.65 & 78.55 & 125.08 & 4.14  & 0.1 \\
			1600  & 78.65 & 78.65 & 151.48 & 2.57  & 0.1 \\
			1800  & 78.7  & 78.7  & 108.57 & 1.72  & 0 \\ 
			\hline
		\end{tabular}
	\end{center}
	\label{tab:14node_cap_data}
\end{table}

\subsection{Design for Maximum Reliability (Topological Expansion)}

\begin{figure}[!hbt]
	\centering
	\includegraphics[width=0.6\textwidth]{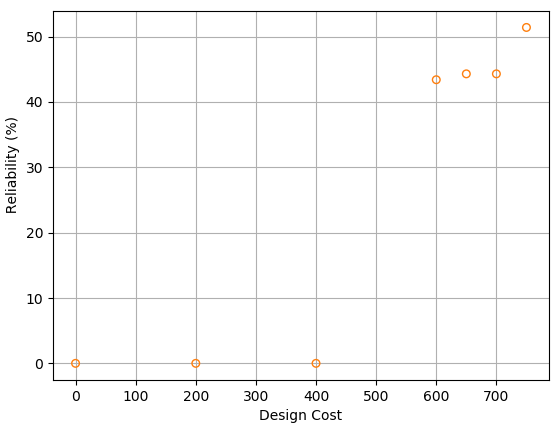}
	\caption{Pareto frontier for optimal topological and capacity design of the 3-node network obtained from the full MILP formulation.}
	\label{fig:3node_top}
\end{figure}

We apply the MILP formulation to the 3-node power network including the use of the topological design variables $v$. We recall that these design variables determine if a particular edge is added to the system. In other words, this more complex formulation chooses an optimal design configuration of edges and capacities where it enforces a fixed upfront cost for the use of each edge. A total of 7 Pareto solutions were computed by varying the budget $\epsilon$ from 0 to 750 and the same samples mentioned above are used. These solutions are presented in Figure \ref{fig:3node_top}. A nonzero $R$ index is not obtained until a budget of 600 is employed since at least 6 edges are required to allow the network to function and the capital cost of each line is 100. After this, capacity increases help improve the network until the $R$ index is maximized by adding all the edges and adding extra capacity resulting in the same best possible optimal design considered shown in Figure \ref{fig:3node_cap_combined}. 

\begin{figure}[!hbt]
	\centering
	\begin{subfigure}[t]{.4\textwidth}
		\centering
		\includegraphics[width=\textwidth]{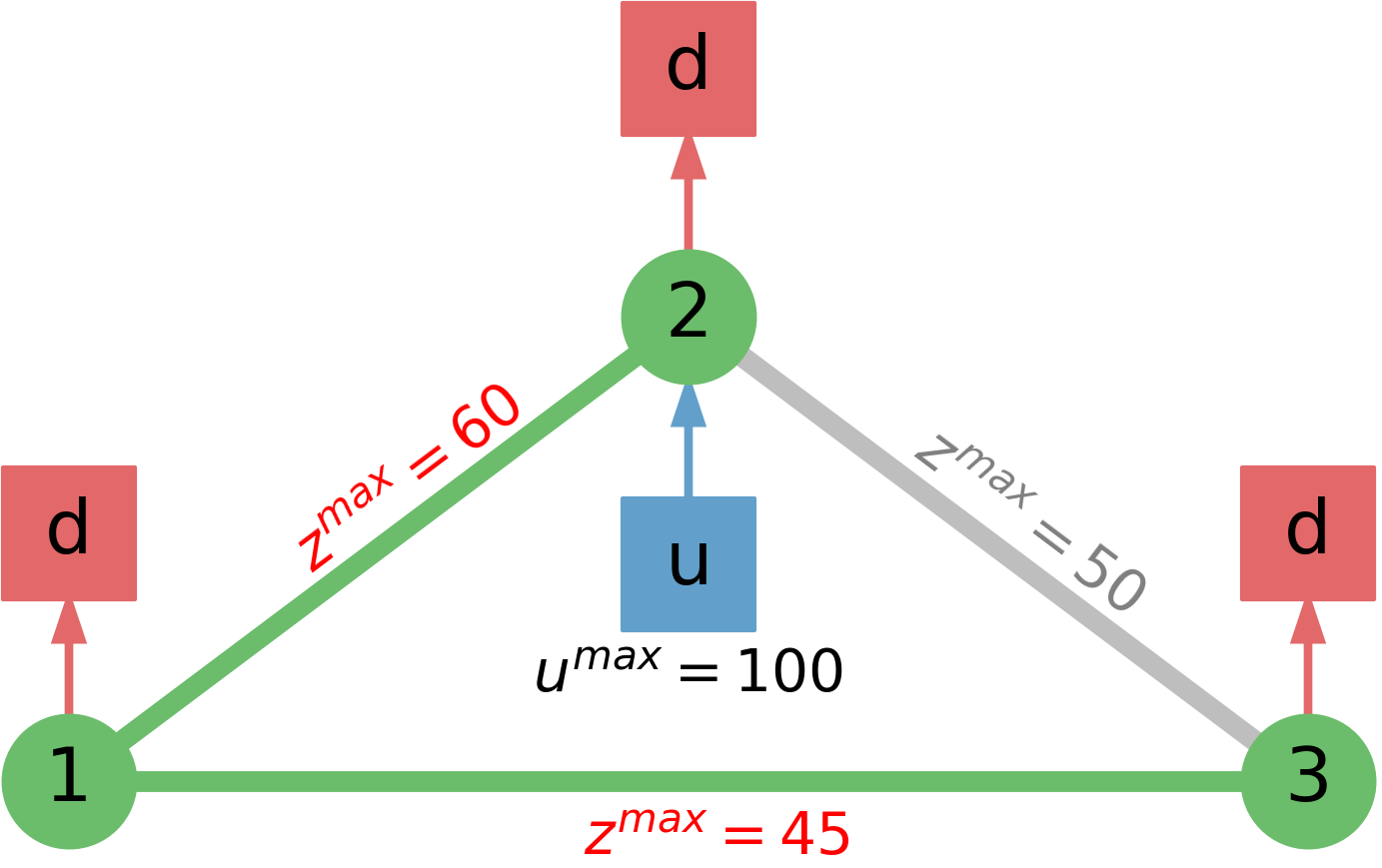}
	\end{subfigure}
	\caption{Schematic of the 3-node network corresponding to the Pareto pair shown in Figure \ref{fig:3node_top} with a design cost of 650.}
	\label{fig:3node_top_combined}
\end{figure}

The optimal design for a budget of 650 is depicted in Figure \ref{fig:3node_top_combined} (left) and this is compared against the design with maximum budget (right). The edges that are switched off are plotted in gray. Here, we observe that the trade-off design is able to effectively operate relative to the sample set without one of the relay edges with the addition of some capacity. The design  of maximum budget employs all of the edges with enough capacity to operate if any relay edge fails (making it the most robust design). 

\begin{figure}[!hbt]
	\centering
	\includegraphics[width=0.6\textwidth]{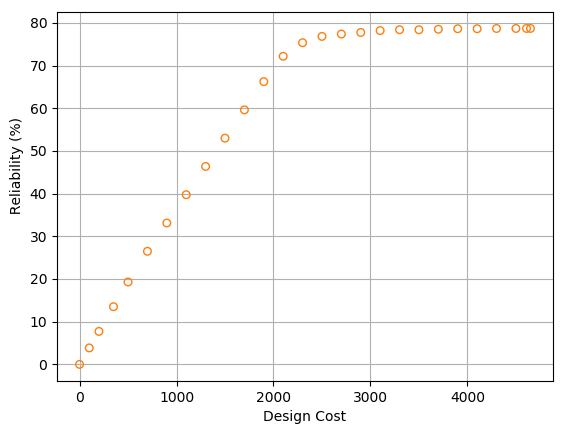}
	\caption{Pareto frontier for optimal topological and capacity design for IEEE 14-node network  (using continuous relaxation).}
	\label{fig:14node_top}
\end{figure}

We also considered topology and capacity design decisions for the IEEE 14-node power network. A total of 27 Pareto pairs were obtained by varying the budget $\epsilon$ from 0 to 4700. The solutions are shown in Figure \ref{fig:14node_top}. The Pareto frontier was obtained using the continuous relaxation. An average of 56 seconds were needed to compute the frontier. The reduced solution times allow us to explore more designs and to explore the reliability limits of the system. In other words, even if the relaxation cannot be guaranteed to provide an exact solution, it captures general behavior and thus can be used as a exploratory tool. 

\begin{figure}[!hbt]
	\centering
	\includegraphics[width=0.6\textwidth]{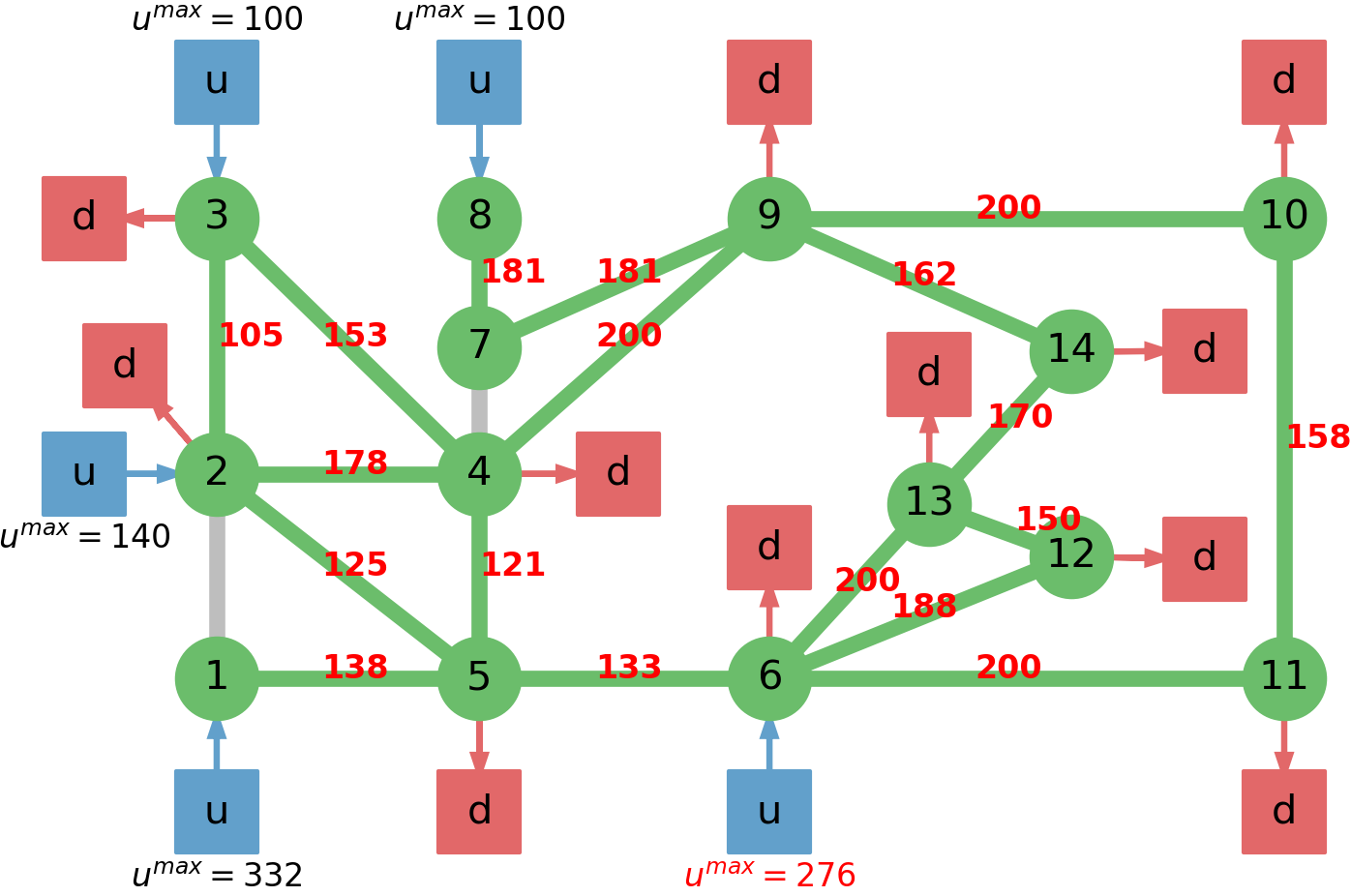}
	\caption{Schematic of the IEEE 14-node network corresponding to Pareto pair shown in Figure \ref{fig:14node_top} with a cost of 2500.}
	\label{fig:14node_top_design}
\end{figure}

The Pareto pair with a  budget of 2500 is shown in Figure \ref{fig:14node_top_design}. The modified capacities are highlighted in red and the edges not in use are colored gray. Here we observe that reliable performance can be obtained simply by adding capacity to the right supplier and most of the edges, except the edges connecting relay nodes 1 to 2 and 4 to 7. Interestingly, this analysis shows that these edges do not impact reliability and can be eliminated (this elimination is not obvious). This highlights that the use of systematic reliability analysis techniques in being can help determine what components are truly needed and avoids over-engineering. 

\section{Conclusions} \label{sec:conclusion}
We propose stochastic programming formulations to compute the reliability of complex systems. Specifically, the proposed reliability measure uses a graph representation of the system and aims to identify the probability that sink nodes are reachable by source nodes. This measure can be computed by solving a network flow problem, can be easily extended to incorporate constraints, and can be easily embedded in design formulations.   We also show that the reliability measure can be computed by solving a stochastic mixed-integer program and that a continuous relaxation of this problem provides high-quality solutions. Case studies are provided to demonstrate the developments. 

\section*{Acknowledgments}
This work was supported by the U.S. Department of Energy under grant DE-SC0014114.

\appendix
\section{Quality of Relaxation for Simple Setting}
\begin{thm}\label{thm:relax_integral} The relaxation of problem \eqref{eq:single_feasibility} is exact. 
\end{thm}

\begin{proof} We express problem \eqref{eq:single_feasibility} (denoted as $P$) in vector form as:
	\begin{equation*}
		\begin{aligned}
		\psi(A,\xi)& =&\max_{y, z} &&&  (1-y) \\
		&&\text{s.t.} &&& A(\xi) z = \hat{d} \cdot (1-y) \\
		&&&&& z \geq 0 \\
		&&&&& y \in \{0, 1\}
		\end{aligned}
	\end{equation*}
	where we define $\hat{d} := -d$ to express the constraints in a standard linear form. The relaxed problem (denoted as $\bar{P}$) is obtained by replacing $y\in\{0,1\}$ with $\bar{y} \in [0, 1]$. We denote optimal solutions for $\bar{P}$ and $P$ as $\bar{y}^*$ and $y^*$, respectively. We show that $\bar{P}$ delivers optimal solutions for $P$ by analyzing three possible cases. First consider the case in which $\hat{d} \in \mathcal{R}(A(\xi))$ (where $\mathcal{R}(\cdot)$ denotes the range of the input matrix) and there exists a nontrivial flow solution ($z^*_j > 0$ for some $j$) such that $A(\xi)z^* = \hat{d}$. This  implies that all values of $y$ are feasible since $(1-y)\hat{d} \in \mathcal{R}(A(\xi)), \ \forall y \in \mathbb{R}$. Thus, $y^* = \bar{y}^* = 0$ must be optimal solutions (yielding the largest possible objective), since any other feasible value of $y$ would have a lower objective. The second case is that in which $\hat{d} \in \mathcal{R}(A(\xi))$ and there does not exist a nontrivial solution ($z^*_j > 0$ for some $j$) such that $A(\xi)z^* = \hat{d}$. In this case the only feasible solution is the trivial solution $z^* = 0$ and  thus $y^* = \bar{y}^* = 1$. 
The third and final case corresponds to $\hat{d} \notin \mathcal{R}(A(\xi))$; it follows that $(1-y)\hat{d} \in \mathcal{R}(A(\xi))$ if and only if $y=1$ since any scalar multiple of $\hat{d}$ will lie outside of $\mathcal{R}(A(\xi))$ except for the trivial case that $\hat{d} = 0$. Thus, the only feasible (and therefore optimal) solution to both problems is $y^* = \bar{y}^* = 1$. 
\end{proof}

\bibliography{references}

\begin{thebibliography}{10}

\bibitem{kim2013network}
Youngsuk Kim and Won-Hee Kang.
\newblock Network reliability analysis of complex systems using a
  non-simulation-based method.
\newblock {\em Reliability engineering \& system safety}, 110:80--88, 2013.

\bibitem{yan2012survey}
Ye~Yan, Yi~Qian, Hamid Sharif, and David Tipper.
\newblock A survey on cyber security for smart grid communications.
\newblock {\em IEEE Communications Surveys \& Tutorials}, 14(4):998--1010,
  2012.

\bibitem{ogunnaike2009random}
Babatunde~A Ogunnaike.
\newblock {\em Random phenomena: fundamentals of probability and statistics for
  engineers}.
\newblock CRC Press, 2009.

\bibitem{thomaidis1994integration}
TV~Thomaidis and EN~Pistikopoulos.
\newblock Integration of flexibility, reliability and maintenance in process
  synthesis and design.
\newblock {\em Computers \& chemical engineering}, 18:S259--S263, 1994.

\bibitem{ye2018mixed}
Yixin Ye, Ignacio~E Grossmann, and Jose~M Pinto.
\newblock Mixed-integer nonlinear programming models for optimal design of
  reliable chemical plants.
\newblock {\em Computers \& Chemical Engineering}, 116:3--16, 2018.

\bibitem{bistouni2014analyzing}
Fathollah Bistouni and Mohsen Jahanshahi.
\newblock Analyzing the reliability of shuffle-exchange networks using
  reliability block diagrams.
\newblock {\em Reliability Engineering \& System Safety}, 132:97--106, 2014.

\bibitem{li2013reliability}
Wenyuan Li et~al.
\newblock {\em Reliability assessment of electric power systems using Monte
  Carlo methods}.
\newblock Springer Science \& Business Media, 2013.

\bibitem{luedtke2008sample}
James Luedtke and Shabbir Ahmed.
\newblock A sample approximation approach for optimization with probabilistic
  constraints.
\newblock {\em SIAM Journal on Optimization}, 19(2):674--699, 2008.

\bibitem{birge2011introduction}
John~R Birge and Francois Louveaux.
\newblock {\em Introduction to stochastic programming}.
\newblock Springer Science \& Business Media, 2011.

\bibitem{straub1993design}
David~A Straub and Ignacio~E Grossmann.
\newblock Design optimization of stochastic flexibility.
\newblock {\em Computers \& Chemical Engineering}, 17(4):339--354, 1993.

\bibitem{bansal1998flexibility}
V~Bansal, JD~Perkins, and EN~Pistikopoulos.
\newblock Flexibility analysis and design of dynamic processes with stochastic
  parameters.
\newblock {\em Computers \& chemical engineering}, 22:S817--S820, 1998.

\bibitem{swaney1985index}
Ross~Edward Swaney and Ignacio~E Grossmann.
\newblock An index for operational flexibility in chemical process design. part
  i: Formulation and theory.
\newblock {\em AIChE Journal}, 31(4):621--630, 1985.

\bibitem{pulsipher2018mixed}
Joshua~L Pulsipher and Victor~M Zavala.
\newblock A mixed-integer conic programming formulation for computing the
  flexibility index under multivariate gaussian uncertainty.
\newblock {\em Computers \& Chemical Engineering}, 119:302--308, 2018.

\bibitem{kim2015guide}
Sujin Kim, Raghu Pasupathy, and Shane~G Henderson.
\newblock A guide to sample average approximation.
\newblock In {\em Handbook of simulation optimization}, pages 207--243.
  Springer, 2015.

\bibitem{hsu1947complete}
Pao-Lu Hsu and Herbert Robbins.
\newblock Complete convergence and the law of large numbers.
\newblock {\em Proceedings of the National Academy of Sciences of the United
  States of America}, 33(2):25, 1947.

\bibitem{kleywegt2002sample}
Anton~J Kleywegt, Alexander Shapiro, and Tito Homem-de Mello.
\newblock The sample average approximation method for stochastic discrete
  optimization.
\newblock {\em SIAM Journal on Optimization}, 12(2):479--502, 2002.

\bibitem{pulsipher2019scalable}
Joshua~L Pulsipher and Victor~M Zavala.
\newblock A scalable stochastic programming approach for the design of flexible
  systems.
\newblock {\em Computers \& Chemical Engineering}, 2019.

\bibitem{dunning2017jump}
Iain Dunning, Joey Huchette, and Miles Lubin.
\newblock Jump: A modeling language for mathematical optimization.
\newblock {\em SIAM Review}, 59(2):295--320, 2017.

\bibitem{zimmerman2010matpower}
Ray~Daniel Zimmerman, Carlos~Edmundo Murillo-S{\'a}nchez, and Robert~John
  Thomas.
\newblock Matpower: Steady-state operations, planning, and analysis tools for
  power systems research and education.
\newblock {\em IEEE Transactions on power systems}, 26(1):12--19, 2010.

\end{thebibliography}
\end{document}